\documentclass[10pt,twoside]{article}
\usepackage{graphicx}
\usepackage{amsmath}
\usepackage{Latex-document}

\newcommand{\goto}{\rightarrow}
\newcommand{\eps}{\epsilon}
\newcommand{\bitem}{\begin{itemize}}
\newcommand{\eitem}{\end{itemize}}

\newcommand{\bR}{{\bf R}}

\newcommand{\cF}{{\cal F}}

\newcommand{\beq}{\begin{equation}}
\newcommand{\eeq}{\end{equation}}

\usepackage{amssymb}
\def\volumeno{I}
\title{\bf  Emerging Applications of \vskip -2mm
Geometric Multiscale Analysis\vskip 6mm}
\author{David L. Donoho\thanks{Department of Statistics,
Stanford University, Stanford, CA 94305, USA.  E-mail:
donoho@stanford.edu}\vspace*{-0.5cm}}
\date{\vspace{-8mm}}

\begin{document}

\maketitle

\thispagestyle{first} \setcounter{page}{209}
\newcommand{\etalchar}[1]{$^{#1}$}
\begin{abstract}\vskip 3mm

Classical multiscale analysis based on wavelets has a number of
successful applications, e.g. in data compression, fast
algorithms, and noise removal. Wavelets, however, are adapted to
point singularities, and many phenomena in several variables
exhibit intermediate-dimensional singularities, such as edges,
filaments, and sheets. This suggests that in higher dimensions,
wavelets ought to be replaced in certain applications by
multiscale analysis adapted to intermediate-dimensional
singularities,

My lecture described various initial
attempts in this direction. In particular,
I discussed two approaches to
geometric multiscale analysis originally arising in the
work of Harmonic Analysts Hart Smith and Peter Jones
(and others): (a)
a directional wavelet transform based on
parabolic dilations; and (b) analysis
via anistropic strips. Perhaps surprisingly,
these tools have
potential applications
in data compression, inverse problems,
noise removal, and signal detection;
applied mathematicians, statisticians,
and engineers are eagerly pursuing these leads.

{\bf Note:} Owing to space constraints, the article is a severely
compressed version of the talk. An extended version of this article,
with figures used in the presentation, is available online at:

\centerline{{\it
http://www-stat.stanford.edu/$\sim$donoho/Lectures/ICM2002}}

\vskip 4.5mm


\end{abstract}

\vskip 12mm

\section{Prologue}

\vskip-5mm \hspace{5mm}

Since the last ICM, we have lost three great mathematical
scientists of the twentieth century:
Alberto Pedro Calder\'on (1922-1999), John Wilder Tukey (1915-2000)
and Claude Elwood Shannon (1916-2001).
%
Although these three are not typically spoken
of as a group, I find it fitting to mention these three
together because  each of these
figures symbolizes for me one aspect of the {\it unreasonable effectiveness
of harmonic analysis}.

Indeed we are all aware of the birth of harmonic analysis
in the nineteenth century
as a tool for understanding of the equations of
mathematical physics, but it is striking how
the original tools of harmonic analysis have frequently
(a) changed, and (b) been applied in ways the inventors
could not have anticipated. Thus, (a) harmonic analysis
no longer means `Fourier Analysis' exclusively,
because wavelet and other forms of decompositions
have been invented by modern harmonic analysts
(such as Calder\'on); and (b) harmonic analysis
finds extensive application outside of mathematical
physics, as a central infrastructural element of
the modern information society, because of the ubiquitous
applications of the fast Fourier transform (after Tukey)
and Fourier transform coding (after Shannon).

There is a paradox here, because harmonic analysts
are for the most part not seeking applications,
or at any rate, what they regard as possible applications
seem not to be the large-scale applications that actually
result. Hence the impact achieved by harmonic analysis
has often not been the intended one
After meditating for a while
on what seems to be the `unreasonable' effectiveness of harmonic analysis,
I have identified what seems to me a chain of argumentation
that renders the `unreasonable' at least `plausible'.
The chain has two propositions:
\bitem
\item {\it Information has its own architecture.}
Each data source, whether imagery, sound, text,
has an inner architecture which we should attempt to
discover and exploit for applications
such as noise removal, signal recovery,
data compression, and fast computation.
\item {\it Harmonic Analysis is about inventing and
exploring architectures for information.} Harmonic analysts
have always created new architectures for decomposition,
rearrangement and reconstruction of operators and functions.
\eitem
In short, {\it the inventory of architectures
created by harmonic analysis amounts to an {\bf intellectual
patrimony} which modern scientists and engineers can fruitfully draw upon
for inspiration as they pursue applications}.
Although there is no necessary connection
between the architectures that harmonic analysts are studying and
the architectures that information requires, it is important
that we have many examples of useful architectures available,
and harmonic analysis provides many of these.  Occasionally, the
architectures already inventoried by harmonic analysts will be exactly
the right ones needed for specific applications.

I stress that the `externally professed goals'
of harmonic analysis in recent decades have always been theorems,
e.g. about the almost everywhere convergence of Fourier Series,
the boundedness of Bochner-Riesz summation operators, or the
boundedness of the Cauchy integral on chord-arc curves.
These externally professed goals have, as far as I know,
very little to do with applications where
harmonic analysis has had wide scale impact.
Nevertheless, some harmonic analysts
are aware of the architectural element in what they
do, and value it highly. As R.R. Coifman has pointed
out to me in private communication:

\begin{quotation}
``The objective of Zygmund, Calder\'on and their school was not
the establishment of new theorems by any means possible. It was
often to take known results that seemed like magic --- e.g.
because of the way they used complex variables methods --- and
tear them apart, finding the underlying structures and their inner
interactions that made it absolutely clear what was going on.  The
test of understanding was measured by the ability to prove an
estimate."
\end{quotation}

In short, the goal was to find the right architecture,
not merely to find the right estimate.

\section{Overview}

\vskip-5mm \hspace{5mm}

In my lecture, I was able  to discuss the possibility
that a coherent subject of {\it Geometric} Multiscale Analysis (GMA)
can be developed -- a subject spanning both mathematics
and a wide range of applications. It is at this point unclear
what the boundaries of the subject will be,
but perhaps the speculative nature of what I had to say
will excite the interest of some readers.  I found it useful
to organize the presentation around the Calder\'on reproducing formula,
which gave us the continuous wavelet transform, but also can be adapted to
give us other multiscale transforms with interesting geometric aspects.
The several different information architectures I described
give an intuitive understanding of what GMA might consist
of.  In the article below, I will review some of the achievements of
classical $1$-dimensional multiscale analysis (wavelet analysis)
starting in the 1980's,
both the mathematical achievements and the extensive applications;
then I will as a warm-up discuss reasons that we
need alternatives to $1$-dimensional multiscale analysis
and its straightforward $d$-dimensional extensions,
and some ideas such as ridgelets, that point in
the expected directions.
In my lecture, I was able to discuss two harmonic analysis
results of the 1990's -- Hart Smith's ``Hardy space
for FIO's'' and Peter Jones' ``Travelling Salesman'' theorem.
Both results concern the higher-dimensional setting,
where it becomes possible to bring in geometric ideas.
I suggested that, in higher dimensions, there are
interesting, nontrivial, nonclassical, geometric multiscale architectures,
with applications paralleling the one-dimensional case. I was
able to sketch some developing applications of these post-classical
architectures.
If these applications can be developed as extensively as has been done for
classical multiscale analysis, the impacts may be large indeed.
In this article, I really have space only to mention topics growing out of
my discussion of Hart Smith's paper. For an extended version of the article,
covering the talk more fully, see \cite{DWTURL}.

{\it Note:} Below we make a distinction between stylized applications
(idealized applications in mathematical models) and actual
applications (specific contributions to scientific discourse and
technological progress); we always
describe the two in separate subsections.

%

\section{Classical multiscale analysis}

\vskip-5mm \hspace{5mm}

An efficient way to introduce  classical multiscale analysis is to
start from Calder\'on's reproducing formula, or as commonly called
today, the {\it Continuous Wavelet Transform}. We suppose we have
a real-valued function $f:\bR \mapsto \bR$ which we want to
decompose into contributions from various scales and locations. We
take with a {\bf wavelet}, an oscillatory real-valued function
$\psi(t)$ satisfying the Calder\'on admissibility condition
imposed on the Fourier transform $\hat{\psi}$ as $\int_0^\infty
|\hat{\psi}(\xi t)|^2 \frac{dt}{t} = 2\pi$, $\forall \xi \neq 0$.
We translate and dilate according to $(\psi_{a,b})(t) = \psi
((t-b)/a) / \sqrt{a}$. We perform {\bf Wavelet Analysis} by
`hitting' the function against all the different wavelets,
obtaining $W_f(a,b) = \langle  \psi_{a,b},  f \rangle$; $W_f$ is
called the Continuous Wavelet Transform (CWT). The CWT contains
all the information necessary to reconstruct $f$, so we can
perform {\bf Wavelet Synthesis} by integrating overall all scales
and locations, summing up wavelets with appropriate coefficients.
\[
f(t) = \int W_f(a,b) \psi_{a,b}(t) \mu(da db).
\]
Here $\mu(da db)$ is the appropriate reference measure,
in this case $\frac{db}{a} \frac{da}{a}$. The `tightness' of the
wavelet transform
as a characterisation of the properties of $f$ is expressed by the
Parseval-type relation $\int W_f(a,b)^2  \mu(da \, db) = \int f(t)^2 dt$.
See also \cite{TenLectures,FJW,Meyer}.

\subsection{Mathematical results}

\vskip-5mm \hspace{5mm}

The CWT maps $f$ into a time-scale plane;
by measuring properties of this time-scale portrait
we can obtain
norms on functions which lead to interesting
theories of functional spaces
and their properties; there are two broad scales of such
spaces we can describe. To define the
{\bf Besov}  $B_{p,q}^\sigma$ spaces we integrate over locations first,
and then over scales
\[
\left(  \int \left( \int (|W(a,b)|a^{-s})^p \frac{db}{a}
\right)^{q/p} \frac{da}{a} \right)^{1/p}.
\]
To define the {\bf Triebel-Lizorkin} $F_{p,q}^\sigma$ spaces
we integrate over scales first and then over locations
\[ \left(  \int \left( \int
(|W(a,b)|a^{-s})^q
\frac{da}{a^{1+q/p}}
\right)^{p/q}
  db \right)^{1/p}.
\]
Here $s = \sigma - (1/p-1/2)$, and we adopt a convention here
and below of {\it ignoring the low frequencies} so that actually
these formulas are only correct for functions which are built from
frequencies $|\xi| > \lambda_0$; the correct general formulas would
require an extra term for the low frequencies which
will confuse the novice and be tedious for experts.
Also for certain combinations of parameters $p,q=1,\infty$
for example, changes ought to be made, based on maximal
functions, BMO norms, etc., but in this expository
work we gloss over such issues.

Each of these norms asks that the wavelet transform decay as
we go to finer scales, and so controls the oscillations of the functions.
Intuition about some of these
spaces comes by thinking of a wavelet coefficient
as something akin to a difference operator such as
$f(b+ a) - 2 f(b) + f(b - a)$; the various
norms on the continuous wavelet coefficients
measure explicitly the finite differences
and implicitly the derivatives
of the analyzed functions.
The distinctions between spaces come in the subtle aspects of choice of
order of integrating in scale and in location
and in choice of $p$ and $q$. We get the following sequence
of relations between the spaces defined by the $F$ and $B$ scales
and classical spaces:
\bitem
\item { $L^p$ :} $L^p \sim F_{p,2}^0$, $1 < p < \infty$.
\item { $H^p$ :} $H^p \sim F_{p,2}^0$, $0 < p \leq 1$.
\item {\bf Sobolev:} $W_p^m \sim F_{p,2}^m$, $1 < p < \infty$.
\item {\bf H\"older:} $C^{\alpha} \sim B_{\infty,\infty}^{\alpha}$
\eitem
There are also equivalences with non-classical, but very interesting, spaces,
such as the Bump Algebra $B_{1,1}^1$, and almost-equivalences
to some other fundamental spaces, such as  $BV(\bR)$.
The full story about such equivalence is told very well
in \cite{Meyer,FJW}.

An important
structural fact about
these spaces is that they admit molecular
decompositions; we can define molecules
as functions obeying certain size, smoothness
and vanishing moment conditions, which are localized
near an interval of some scale and location,
and then show that, although elements of these spaces
are defined by norms on the continuum domain,
functions belong to these spaces if and only if
they can be written as superpositions
$ f(x) = \sum_Q A_Q m_Q(x)$ where $m_Q$ are molecules
and the $A_Q$ are scalar coefficients, and
where the coefficient sequence $(A_Q)_Q$ obeys certain
norm constraints. Results of this kind first
emerged in the 1970's;
a canonical way to get such results uses the CWT \cite{FJW}.
Consider the dyadic cells
\[
Q = \{ (a,b):  2^{-j} > a \geq 2^{-(j+1)}, k/2^j \leq b < (k+1)/2^j \},
\]
note that they obey $\mu(Q) \approx 1$; they are ``unit cells'
for the reference measure. It turns out that
the behavior of $W(a,b)$
at various points within such a cell $Q$
stays roughly comparable \cite{TenLectures,FJW},
and that the $\psi_{a,b}$ all behave similarly as well.
As a result, the
integral decomposition offered by the Calder\'on reproducing formula
can sensibly be discretized into terms arising from different cells.
\begin{eqnarray*}
    f(x) &=& \int W(a,b) \psi_{a,b}(x) \mu(da \,db) \\
         &=& \sum_Q \int_Q W(a,b) \psi_{a,b}(x) \mu(da \,db)\\
         &=& \sum_Q M_Q(x), \qquad M_Q = \int_{Q} W(a,b) \psi_{a,b}
\mu(da \,db) \\
         &=& \sum_Q A_Q m_Q(x), \qquad A_Q = \|W(\cdot,\cdot) \|_{L^2(Q)}
\end{eqnarray*}
Now roughly speaking, each $m_Q$ is a mixture of wavelets
at about the same location and scale, and so is something
like a wavelet, a coherent oscillatory waveform of a certain
location and scale.
This type of discretization of the Calder\'on reproducing formula
has been practiced since the 1970's,
for example by Calder\'on,
and by Coifman and Weiss \cite{CoifmanAtomic,CoifmanWeiss},
who introduced the terms {\it molecular decomposition} (and atomic
decomposition) for discrete series of terms
localized near a certain scale and location.
Hence, the $m_Q$ may be called
molecules and the $A_Q$ represent
the contributions of various molecules.
The spaces $F_{p,q}^\sigma$
and $B_{p,q}^\sigma$ can then be characterized by the decomposition
$ f(x) = \sum_Q A_Q m_Q(x)$: we can define sequence-space
norms $f_{p,q}^\sigma$ as in (\ref{triebseqdef}) below for which
\[
    \| f \|_{F_{p,q}^\sigma} \sim \| (A_Q)_Q \|_{f_{p,q}^\sigma},
\]
and similarly for Besov sequence norms $b_{p,q}^\sigma$.
This gives a clear understanding  of the structure of $f$
in terms of the distribution of the number and size of oscillations across
scales.

While the molecular decomposition
is very insightful and useful for proving
structure theorems about functional spaces,
it has two drawbacks which severely restrict practical applications.
First, the $A_Q$ are nonlinear functionals of the underlying object $f$;
secondly, the $m_Q$ are variable objects which depend on $f$.
As a result, practical applications of the sum $\sum_Q A_Q m_Q$
are not as straightforward
as one might like. Starting in the early 1980's,
it was found that a much simpler and more practical
decomposition was possible;
in fact with appropriate choice of generating
wavelet -- different than usually made in the CWT --
one could have an {\it orthonormal wavelet basis}
\cite{Meyer,TenLectures}.
\begin{equation} \label{orthodef}
      f = \sum_{j,k} W(2^{-j},k/2^j)  \psi_{2^{-j},k/2^j} = \sum_{j,k}
\alpha_{j,k} \psi_{j,k}.
\end{equation}
Essentially, instead of integrating over dyadic cells $Q$,
it is necessary only to sample once per cell! Several crucial advantages
in applications flow from the fact that the coefficients $\alpha_{j,k}$
are linear in $f$ and the $\psi_{j,k}$
are fixed and known.

A theoretical advantage flows from the fact that
the same norm equivalence that was
available for the amplitudes $(A_Q)$ in
the molecular decomposition also applies
for the wavelet coefficients:
\begin{eqnarray}
      \| f \|_{B_{p,q}^\sigma} &\sim&  \| \alpha \|_{b_{p,q}^\sigma}
\equiv \left( \sum_j ( \sum_k
|\alpha_{j,k}|^p )^{q/p} 2^{jsq}
\right)^{1/q}  \label{triebseqdef},\\
      \| f \|_{F_{p,q}^\sigma} &\sim&  \| \alpha \|_{f_{p,q}^\sigma}
\equiv  \left( \int ( \sum_j  |\alpha_{j,k}|^q 2^{jsq}
\chi_{j,k}(t) )^{p/q}  \label{besovseqdef} \right)^{1/q}.
\end{eqnarray}
This implies that the wavelets $\psi_{j,k}$ make an unconditional
basis for appropriate spaces in the Besov and Triebel scales. This
can be seen from the fact that the norm involves only
$|\alpha_{j,k}|$; a fact which is quite different from the case
with Fourier analysis. Unconditionality implies that the balls $\{
f: \| f \|_{B_{p,q}^\sigma} \leq A \}$ are closely inscribed by
and circumscribed by balls $\{ f : \| \alpha \|_{b_{p,q}^\sigma}
\leq A' \}$ which are quite simple geometric objects, solid and
orthosymmetric with respect to the wavelets as `principal axes'.
This solid orthosymmetry is of central significance for the
optimality of wavelets for many of the stylized applications
mentioned below; compare \cite{UBasis,Bits,DJMinimax}.

Our last chapter in the mathematical development
of classical multiscale methods
concerns the connection
between Besov spaces and approximation spaces.
In the late 1960's, Jaak
Peetre observed that the space
$B_{1/\sigma,1/\sigma}^\sigma$, $\sigma > 1$ was very special.
It served as the {\it Approximation Space}
for approximation in $L^\infty$ norm by
free knot splines, i.e. as the set of functions
approximable at rate $n^{-\sigma}$ by splines with $n$ free knots.
In the 1980's a more general
picture emerged,  through work of e.g. Brudnyi,
  DeVore, Popov and Peller \cite{DeVorePopov}: that the space
$B_{\tau,\tau}^{\sigma}$ served as the
approximation space of many nonlinear approximation schemes
(e.g. rational functions), under $L^p$ approximation error,
where $1/\tau = \sigma + 1/p$.  This says that although $\tau < 1$
at first seems unnatural (because graduate mathematical training emphasizes
convex spaces) these nonconvex spaces are fundamental.
The key structural fact is that those spaces are equivalent,
up to renorming, to the set of functions whose wavelet
coefficients belong to an $\ell^\tau$ ball, $\tau < 1$.
Hence, membership of wavelet coefficients in an $\ell^\tau$ ball
for small $\tau$ becomes of substantial interest. The intuitive
appeal for considering $\ell^\tau$ balls is clear by considering
the closely-related weak-$\ell^\tau$ balls;  they can be defined
as the constants $C$ in relations of the form
\[
    \mu \{ (a,b): |W(a,b)| > \eps \} \leq C \eps^{-1/\tau}, \quad \eps > 0,
\]
or
\[
    \# \{ (j,k): |\alpha_{j,k}| > \eps \} \leq C \eps^{-1/\tau}, \quad \eps > 0.
\]
They are visibly measures of the sparsity in the time-scale plane,
and hence sparsity of that plane controls the asymptotic behavior of
numerous nonlinear
approximation schemes.

\subsection{Stylized applications}

\vskip-5mm \hspace{5mm}

We now mention some stylized applications of classical multiscale
thinking, i.e. applications in a model world where we can
prove theorems in the model setting.

\subsubsection{Nonlinear approximation}

\vskip-5mm \hspace{5mm}

Since the work of D.J. Newman in the 1960's it
was understood that approximation by rational
functions could be dramatically better than
approximation by polynomials; for example
the absolute value function $|t|$ on the interval
$[-1,1]$ can be approximated at an exponential
rate in $n$ by rational functions with numerator and
denominator of degree $n$,
while it can be approximated only at an algebraic rate $n^{-1}$
by polynomials of degree $n$. While this suggests the
power of rational approximation, it must also
be noted that rational approximation is a
highly nonlinear and computationally complex
process.

On the other hand, from
the facts (a) that wavelets provide an unconditional basis for
Besov spaces, and (b) that certain Besov spaces are approximation
spaces for rational approximation,
we see that wavelets give an effective
algorithm for the same problems
where rational functions would be useful.
Indeed, based on work by
DeVore, Popov, Jawerth, Lucier
we know that if we consider the class of
functions approximable at rate
$\approx n^{-\tau}$ by rational approximation,
these same functions can be approximated
at the same rate simply by taking a partial reconstruction
based on the  $n$
``biggest'' wavelet coefficients.

In short, from the viewpoint of asymptotic rates of convergence,
thresholding of wavelet coefficients \-- a very weakly nonlinear
approximation scheme \-- is fully as effective as best rational
approximation.  The same assertion can be made comparing
nonlinear approximation by wavelets and by free knot splines.

\subsubsection{Data compression}

\vskip-5mm \hspace{5mm}

Consider the following mathematical idealization of data compression. We
have a function which is an unspecified element of a Besov Ball $ \cF =
\{ f : \| f \|_{B_{p,q}^\sigma} \leq A \}$ and we wish to have a
coder/decoder pair which can approximate any such function to within an
$\eps$-distance in $L^2$ norm by encoding into, and decoding from, a
finite bitstring.

In mathematical terms, we are studying the
Kolmogorov $\eps$-entropy: we wish to achieve $N(\eps,\cF)$, the
minimal number of bits required to represent every $f$ in $\cF$
to within an $L^2$ error $\eps$.
This is known, since Kolmogorov and Tikhomirov, to behave as
\begin{equation} \label{optrate}
      N(\eps, \cF) \asymp \eps^{-1/\sigma}, \eps \goto 0.
\end{equation}
Now, up to renorming,
the ball $\cF$ is isometric to a ball
in sequence space $\Theta = \{ \alpha : \| \alpha \|_{b_{p,q}^\sigma}
\leq A \}$.
Such a ball is a subset of $w\ell^\tau$
for $1/\tau = \sigma + 1/2$
and each element in it
can be approximated in $\ell^2$ error
at a rate $M^{-1/\tau+1/2}$ by
sparse vectors containing only $M$
nonzero coefficients. Here is a simple
coder inspired by this fact. Pick $M(\eps)$ coefficients such that
the $\ell^2$-error of such an approximation is
at most (say) $\eps/2$. The $M(\eps)$ coefficients
achieving this can be
quantized into integer multiples of a base quantum $q$,
according to $a_{j,k} = \lfloor \alpha_{j,k} / q \rfloor$,
with the
quantum chosen so that the quantized vector $\alpha^{(q)}$ defined
by  $ \alpha^{(q)}_{j,k} = q \cdot a_{j,k}$, approximates the original
coefficients to within $\ell^2$ error $\eps/2$.
The resulting integers $a_{j,k}$ represent the function $f$
to within $L^2$ error $\eps$
and their indices can be coded into bit strings,
for a total encoding length of not worse that
$O(\log(\eps^{-1}) M(\eps)) = O(\log(\eps^{-1}) \eps^{-1/\sigma})$.
Hence a very simple algorithm on the wavelet coefficients
gets close to the optimal asymptotics (\ref{optrate})! Underlying
this fact is the geometry of the body $\Theta$; because of its solid
orthosymmetry, it  contains many high-dimensional
hypercubes of ample radius. Such hypercubes are
essentially incompressible.

In fact the $\log(\eps^{-1})$ factor is removable
in a wide range of $\sigma,p,q$.
In many cases, the $\eps$-entropy can be attained,
within a constant factor, by appropriate level-dependent
scalar
quantization of the wavelet coefficients
followed by run-length encoding.
In other work, Cohen et al. have shown that
by using the tree-organization of wavelet coefficients
one can develop algorithms which give the right
order of asymptotic behavior for the across
many smoothness classes; e.g. \cite{WaveCompress}.

In fact more is true.
Suppose we use for Besov ball
simply the ball $ \{ f : \| \alpha(f) \|_{b_{p,q}^\sigma} \leq A \}$
based on wavelet coefficients;
then by transform coding as
in \cite{CountingBits}
we can get efficient
codes with codelength precisely asymptotic equivalence to the
Kolmogorov $\eps$-entropy by levelwise
$\ell^p$-sphere vector quantization of wavelet coefficients.
Underlying this fact,
the representation of the underlying functional
class as an orthosymmetric body
in infinite-dimensional space is very important.

\subsubsection{Statistical estimation}

\vskip-5mm \hspace{5mm}

Consider the following mathematical idealization of nonparametric curve
estimation. We have an unknown function
  $f(t)$ on $[0,1]$  which is an element of a Besov Ball $ \cF = \{ f : \| f
\|_{B_{p,q}^\sigma}
\leq A
\} $
We observe data $Y$ from the white noise model
\[
    Y(dt) = f(t)dt + \eps W(dt),
\]
where the $W(t)$ is a Wiener process and
$\eps$ the noise level,
and we wish to reconstruct $f$ accurately.
We measure risk using
the mean squared error
\[
    R_\eps(f,\hat{f})  = E \| f - \hat{f} \|_2^2.
\]
and evaluate quality by the minimax risk
\[
      \min_{\hat{f}} \max_{f \in \cF} R_\eps(f,\hat{f}) .
\]
Over a  wide range of $\sigma, p, q$,
this minimax risk tends to zero at the rate $(\eps^2)^{2\sigma/(2\sigma+1)}$.

In this setting, some simple
algorithms based on noisy wavelet coefficients $y_{j,k} = \int
\psi_{j,k}(d) Y(dt)$
can be quite effective. In effect, $y_{j,k} = \alpha_{j,k} + \eps z_{j,k}$,
where $z_{j,k}$ is a white Gaussian noise. By simply applying
thresholding
to the noisy wavelet coefficients of $Y$,
\[
\hat{\alpha}_{j,k} = y_{j,k} 1_{\{ |y_{j,k}| > \lambda \eps \}}
\]
at scales $0 \leq j \leq \log_2(\eps^{-2})$
with threshold $\sqrt{2 \log(\eps^{-1})}$ ,
  we obtain a new set of coefficients; using these
we obtained a nonlinear approximation $\hat{f} = \sum_{j,k}
\hat{\alpha}_{j,k} \psi_{j,k}$. The quantitative properties are
surprisingly good; indeed, using again the $w\ell^\tau$ embedding
of the Besov body $b_{p,q}^\sigma$, we have that the
$\ell^2$-error of nonlinear approximation to $\alpha$ using $M$
terms converges at rate $M^{-1/\tau + 1/2}$. Heuristically, the
coefficients surviving thresholding have errors of size $\approx
\eps$, and the object can be approximated by at most $M$ of these
with $\ell^2$ error $\approx M^{-1/\tau + 1/2}$; simple
calculations suggest that the risk of the estimator is then
roughly $\eps^2 \cdot M + M^{-2/\tau + 1}$ where $M$ is the number
of coefficients larger than $\eps$ in amplitude; this is the same
order as the minimax risk $\eps^{2\sigma/(2\sigma+1)}$! (Rigorous
analysis shows that for this simple algorithm, log terms intervene
\cite{DJKPAsymptopia}.) If we are willing to refine the
thresholding in a level-dependent way, we can obtain a risk which
converges to zero at the same rate as the minimax risk as $\eps
\goto 0$, e.g. \cite{DJMinimax}. Moreover, if we are willing to
adopt as our Besov norm the sequence space $b_{p,q}^\sigma$ norm
based on wavelet coefficients, then by applying a sequence of
particular scalar nonlinearities to the noisy wavelet coefficients
(which behave qualitatively like thresholds) we can get precise
asymptotic equivalence to the minimax risk, i.e. precise
asymptotic minimaxity \cite{DJMinimax}. Parallel results can be
obtained with wavelet methods in various inverse problems, where
$f$ is still the estimand, but we observe noisy data on $Kf$
rather than $f$, with $K$ a linear operator, such as convolution
or Radon transform \cite{WVD}. \pagebreak
\subsubsection{Fast computation}

\vskip-5mm \hspace{5mm}

An important theme for scientific computation
is the sparse representation, not of functions, but of operators.
For this purpose a central fact pointed
out by Yves Meyer \cite{Meyer} is that wavelets sparsify large classes of
operators. Let $T$ be a Calderon-Zygmund operator (CZO);
the matrix representation of such operator in the wavelet
basis
\[
    M_{j,k}^{i,l} = \langle \psi_{j,k} , T \psi_{i,l} \rangle ,
\]
then $M$ is sparse -- all its rows and columns have
finite $\ell^p$ norms for each $p > 0$. In short,
such an operator involves interactions between very few
pairs of terms.

For implications of such sparsity, consider
the work of Beylkin, Coifman, and Rokhlin \cite{BCR}.
Suppose $T$ is a CZO,
and let $Comp(\eps,n)$ denote the
number of flops required to compute
an $\eps$-approximation to $P_n T P_n$,
where $P_n$ is an projector onto scales
larger than $1/n$.  In \cite{BCR}
it was shown that, ignoring set-up costs,
\[
    Comp(\eps,n) = O( \log(1/\eps) n )
\]
so that such operators could be applied many
times with cost essentially linear
in problem size, as opposed to the $O(n^2)$ cost
nominally demanded by matrix multiplication.
The algorithm was roughly this: represent the
operator in a wavelet basis, threshold the coefficients,
and keep the
large coefficients in that representation.
A banded matrix results, which can be applied in order $O(n)$
flops.  (The story is a bit more subtle, since the
algorithm as written would suffer an additional
$O(\log(n))$ factor; to remove this, Beylkin,
Coifman, and Rokhlin's nonstandard form must be applied.)

\subsection{Applications}

\vskip-5mm \hspace{5mm}

The possibility of applying wavelets to real problems relies heavily on
the breakthrough made by Daubechies \cite{Daubechies} (building on work
of Mallat \cite{Mallat}) which showed that it was possible to define a
wavelet transform on finite digital signals which had orthogonality and
could be computed in order $n$ flops. Once this algorithm was available,
a whole range of associated fast computations followed. Corresponding to
each of the `stylized applications' just listed, many `real
applications' have been developed over the last decade; the most
prominent are perhaps the use of wavelets as part of the JPEG-2000 data
compression standard, and in a variety of signal compression and
noise-removal problems. For reasons of space, we omit details, referring
the reader instead to \cite{DCHA} and to various wavelet-related
conferences and books.

\section{Need for geometric multiscale analysis}

\vskip-5mm \hspace{5mm}

The many successes of classical multiscale analysis
do not exhaust the opportunities for successful multiscale
analysis. The key point is the slogan we formulated
earlier \-- {\it Information has its own architecture}.
In the Information Era, where  new data sources are proliferating
endlessly, each with its own peculiarities and specific phenomena,
there is a need for expansions uniquely adapted to
each type of data.

In this connection, note that classical wavelet analysis
is uniquely adapted to objects which are smooth
apart from {\it point singularities}.  If a function is
$C^\infty$ except for step discontinuities at a finite
set of points, its continuous wavelet transform will be very
sparse. In consequence, the decreasing rearrangement
of its wavelet coefficients will decay
rapidly, and $n$-term approximations to the object will
converge rapidly in $L^2$ norm. With the right definitions
the story in high dimensions is similar: wavelets
give a sparse representation of point singularities.

On the other hand, for singularities along
lines, planes, curves, or surfaces,
the story is quite different.
For functions in dimension $2$
which are discontinuous along a curve, but otherwise
smooth, the $2$-dimensional CWT will not be sparse.
In fact, the the decreasing rearrangement
of its wavelet coefficients will decay
like $C/N$, and $N$-term approximations to the object will
converge no faster than $O(N^{-1})$ in squared $L^2$ norm.
Similar statements can be made for singularities of dimension
$0 < k < d$ in dimension $d$.
In short, wavelets are excellent for representing
smooth data containing point singularities
but not singularities of intermediate dimensions.

There are many examples of data where singularities
of intermediate dimensions constitute important features.
One example comes from extragalactic astronomy, where
gravitational clustering has caused matter to congregate in
`filaments' and 'sheets' in 3-dimensions.  Another example
comes from image analysis, say of SAR imagery,
where stream beds, ridge lines, roads and
other curvilinear phenomena punctuate
the underlying background texture.
Finally, recently-developed tools for $3D$ imaging
offer volumetric data of phsyical objects
(eg biological organs) where sheetlike
structures are important.

We can summarize our vision for
the future of multiscale analysis
as follows.

{\sl If it is possible to
sparsely analyze objects which are
smooth apart from
intermediate-dimensional singularities,
this may open
{\bf new vistas in mathematical
analysis}, offering
(a) new functional Spaces,
and (b) new representation of
mathematically important
operators.

If, further, it is possible
algorithmize such analysis tools,
this would open new applications involving
(a) data compression; (b) noise removal and
recovery from Ill-posed inverse problems;
(c) feature extraction and pattern recognition;
and (d) fast solution of differential and integral equations.
}

But can we realistically expect to sparsely analyse
such singularities? By considering
Calder\'on-like formulas,
we can develop some understanding.

\subsection{Ridgelet analysis}

\vskip-5mm \hspace{5mm}

We consider first the case of singularities of co-dimension 1. It turns
out that the ridgelet transform is adapted to such singularities.

Starting from an admissible
wavelet $\psi$, define the ridgelet
$\rho_{a,b,\theta}(x) = \psi_{a,b}(u'_\theta x) $,
where $u_\theta$ is a unit vector
pointing in direction $\theta$
and so this is a wavelet in one direction and constant in
orthogonal directions \cite{HarmNet}.
In analogy to the continuous
wavelet transform, define the continuous
ridgelet transform
$R_f(a,b,\theta) = \langle \rho_{a,b,\theta},  f
\rangle$.  There is a synthesis formula
\[
  f(x) = \int R_f(a,b,\theta) \rho_{a,b,\theta}(x) \mu(da \,db \,d\theta)
\]
and a Parseval relation
\[
  \|f \|_2^2 = \int R_f(a,b,\theta)^2 \mu(da \,db \,d\theta)
\]
both valid for an appropriate reference measure $\mu$.
Note the similarity to the Calder\'on formula.

In effect this is an analysis of $f$ into
contributions from `fat planes'; it has been
extensively developed in
Emmanuel Cand\`es'
Stanford thesis (1998) and later publications.
Suppose we use it to analyze a
function $f(x) \in L^2({\bf R}^n)$
which is smooth apart from a singularity across
a hyperplane.
If our function is, say,
$f_{u,a}(x) = 1_{\{u'x > a\}} e^{-\|x\|^2}$,
Cand\`es \cite{CandesThesis}.
showed that the ridgelet transform of $f_{u,a}$ is sparse.
For example, a sampling of the continuous ridgelet transform at dyadic
locations and scales and directions gives a set of coefficients
such that the rearranged ridgelet coefficients decay rapidly.
It even turns out that we can define ``orthonormal ridgelets''
(which are not true ridge functions)
such that the orthonormal ridgelet coefficients  are sparse:
they belong to every $\ell^p$ with $p > 0$ \cite{OrthoRidgelets}.
In short, an appropriate multiscale analysis
(but not wavelet analysis) successfully
compresses singularities of co-dimension
one.

\subsection{\boldmath $k$-plane ridgelet transforms}

\vskip-5mm \hspace{5mm}

We can develop comparable reproducing formulas of co-dimension $k$ in
$R^d$. If $P_k$ denotes orthoprojector onto a $k$-plane in $R^d$, and
$\psi$ an admissible wavelet for $k$-dimensional space, we can define a
$k$-plane Ridgelet: $\rho_{a,b,P_k}(x) = \psi_{a,b}(P_k x) $ and obtain
a $k$-plane ridgelet analysis: $R_f(a,b,P_k) = \langle \rho_{a,b,P_k},
f \rangle$. We also obtain a reproducing formula
\[
  f(x) = \int R_f(a,b,P_k) \rho_{a,b,P_k}(x) \mu(da \,db \,dP_k)
\]
and a Parseval relation  $\|f \|_2^2 = \int R_f(a,b,P_k)^2 \mu(da
\,db \,dP_k)$,
with in both cases $\mu()$ the appropriate reference measure. In short we are
analyzing the object $f$ into
`Fat Lines', `Fat $k$-planes,' $ 1 \leq k \leq n-1$.
Compare \cite{k-plane-ridgelet}.
Unfortunately, all such representations
have drawbacks, since to use them
one must fix in advance the co-dimension
$k$; moreover, very few singularities are globally flat!

\subsection{Wavelet transforms for the full affine group}

\vskip-5mm \hspace{5mm}

A more ambitious approach is to consider wavelets indexed by the general
affine group $GA(n)$; defining $(\psi_{A,b} g)(x) = \psi(Ax+b) \cdot
|A|^{1/2}$. This leads to the wavelet analysis $W_f(A,b) = \langle
\psi_{A,b},  f \rangle$. Taking into account the wide range
of'anisotropic dilations and directional preferences poossible within
such a scheme, we are analyzing $f$ by waveforms which represent a very
wide range of behaviors:
  `Fat Points',  `Fat Line Segments',  `Fat Patches',
and so on.

This exciting concept unfortunately fails.
No matter what wavelet we pick to begin with,
$\int W_f(A,b)^2  \mu(dA db) = +\infty$.
(technically speaking, we cannot get a
  square-integrable representation of the general
affine group; the group is too large)
\cite{Weiss1,Weiss2}.
Moreover,
synthesis fails: $\int W_f(A,b) \psi_{a,b}(t) \mu(dA \, db)$
is not well-defined. Finally, the transform is not sparse
on singularities.

In short, the dream of
using Calder\'on-type formulas to
easily get a decomposition of piecewise
smooth objects into  `Fat Points',
  `Fat Line Segments',
`Fat Surface Patches',
and so on fails.
Success
will require hard work.

\subsection{A cultural lesson}

\vskip-5mm \hspace{5mm}

The failure of soft analysis is not unexpected, and not catastrophic. As
Jerzy Neyman once said: life is complicated, but not uninteresting. As
Lennart Carleson said:
\begin{quotation}
There was a period, in the 1940's and 1950's, when classical analysis
was considered dead and the hope for the future of analysis was considered
to be in the abstract branches, specializing in generalization.  As is now
apparent, the death of classical analysis was greatly exaggerated ...
the reasons
for this ... [include] ... the realization that in many problems complications
cannot be avoided, and that intricate combinatorial arguments rather
than polished
theories are in the center.
\end{quotation}

Our response to the failure of Calder\'on's
formula for the full $Ax+b$ group
was to consider, in the ICM Lecture, two specific strategies
for decomposing multidimensional objects.
In the coming section,
we will consider analysis using a special subset of the
$Ax+b$ group, where a Calder\'on-like formula still applies,
and we can construct a fairly complete analog of the
wavelet transform -- only one which is efficient for
singularities of co-dimension 1.
In the lecture (but not in this article),
we also considered analysis using
a fairly full subset of the $Ax+b$ group, but in a simplified
way, and extracted the results we need by special
strategies (viz. Carleson's ``intricate combinatorial arguments'')
rather than smooth general machinery.  The results
delivered in both approaches seem to indicate the
correctness of the vision articulated above.

\section{Geometric multiscale analysis `with Calder\'on'}

\vskip-5mm \hspace{5mm}

In harmonic analysis since the 1970's
there have been a number of important applications
of decompositions based on
{\it parabolic dilations}
\[
     f_a(x_1,x_2) = f_1(a^{1/2}x_1, a x_2),
\]
so called because they leave invariant the
parabola $x_2 = x_1^2$.
Calder\'on himself used such dilations \cite{CalderonParabolic}
and exhibited a reproducing formula where the scale variable
acted through such dilations. Note that in the above equation
the dilation is always twice as strong in one fixed direction as in the
orthogonal one.

At the same time, decompositions began to be used
based on {\bf directional parabolic dilations} of the form
\[
     f_{a,\theta}(x_1,x_2) = f_a(R_\theta (x_1,x_2)').
\]
Such dilations (essentially) leave invariant curves defined by
quadratic forms  with $\theta$ as one of the principal directions.
For example, Charles Fefferman in effect used decompositions based
on parabolic scaling in his study of Bochner-Riesz summability
cite{Fefferman}. Elias Stein used decompositions exhibiting
parabolic scaling in studying oscillatory integrals in the 1970's
and 1980's \cite{FatStein}. In the 1990's, Jean Bourgain, Hart
Smith, Chris Sogge, and Elias Stein found applications in the
study of oscillatory integrals and Fourier Integral operators.

The principle of parabolic scaling
leads to a meaningful decomposition reminiscent
of the continuous wavelet transform, only with a much more strongly
directional character. This point has been developed in
a recent article of Hart Smith \cite{SmithFIO},
who defined a continuous wavelet transform based
on parabolic scaling, a notion of directional
molecule, showed that FIO's map directional molecules into
directional molecules, and showed that FIO's have a sparse
representation in a discrete decomposition.
For this expository work, we have developed
what seems a conceptually simple, perhaps
novel way of approaching this topic, which we hope will
be accesible to non-experts. Details underlying the exposition
are available from \cite{DWTURL}.

\subsection{Continuous directional multiscale analysis}

\vskip-5mm \hspace{5mm}

We will work exclusively in $\bR^2$, although everything generalizes to
higher dimensions. Consider a family of directional wavelets with three
parameters: scale $a > 0$, location $b \in \bR^2$ and orientation
$\theta \in [0,2\pi)$. The orientation and location parameters are
defined by the obvious rigid motion
\[
     \psi_{a,b,\theta} = \psi_{a,0,0}(R_\theta (x - b))
\]
with $R_\theta$ the 2-by-2 rotation matrix
effecting planar rotation by $\theta$ radians.
At fine scales,
the scale parameter $a$ acts in a slightly nonstandard
fashion based on parabolic
dilation, {\it in the polar Fourier domain}.
We pick a wavelet $\psi_{1,0}$ with $\hat{\psi}$ of compact support
away from $0$, and a bump $\phi_{1,0}$ supported in $[-1,1]$.
Here $\psi_{1,0}$ should obey the usual admissibility
condition and  $\|\phi\|_2 = 1$.
At sufficiently fine scales (say $a < 1/2$) we define the directional
wavelet by
going to polar coordinates $(r, \omega)$ and setting
\[
    \hat{\psi}_{a,0,0}(r,\omega) = \hat{\psi}_{a,0}(r) \cdot
\phi_{a^{1/2},0}(\omega), \quad a < a_0 .
\]
In effect, the scaling is parabolic in the polar variables $r$ and
$\omega$, with $\omega$ being the `thin' variable;
thus in particular the wavelet $\psi_{a,0,0}$ is
not obtainable by affine change-of-vartiables on $\psi_{a',0,0}$ for
$a' \neq a$.
We omit description
of the transform at coarse scales, and so again ignore low frequency
adjustment terms. Note that it is correct to call these wavelets
directional, since they become increasingly needle-like at fine scales.

Equipped with such a family of high-frequency
wavelets, we can define a {\it Directional Wavelet Transform}
\[
     DW(a,b,\theta) = \langle \psi_{a,b,\theta} , f \rangle, \quad
a >0, b \in \bR^2, \theta \in [0,2\pi)
\]
It is easy to see that we have a Calder\'on-like
reproducing formula, valid for high-frequency functions:
\[
     f(x) = \int DW(a,b,\theta) \psi_{a,b,\theta}(x) \mu(da \, db \, d\theta)
\]
and a Parseval formula for high-frequency functions:
\[
     \|f\|_{L^2}^2 = \int DW(a,b,\theta)^2  \mu(da \, db \, d\theta)
\]
in both cases, $\mu$ denotes the reference measure
$\frac{db}{a^{3/2}} \frac{d\theta}{a^{1/2}} \frac{da}{a}$.

Based on this transform, we can define seminorms
reminiscent of Besov and Triebel seminorms in wavelet analysis;
while it is probably a major task to prove that thse
give well-founded spaces, and such work has not
yet been done (for the most part), it still seems useful
to use these as a tool measuring the distribution of a function's
`content' across scale, location and direction. We get a directional
{\bf Besov}-analog $DB_{p,q}^\sigma$: integrating over locations
and orientations first
\[
\left(  \int \left( \int
(|DW(a,b,\theta)|a^{-s})^p  \frac{d\theta}{a^{1/2}} \frac{db}{a^{3/2}}
\right)^{q/p}
\frac{da}{a^2}
\right)^{1/p}
\]
and a {\bf Triebel}-analog $DF_{p,q}^\sigma$ by integrating over scales first
\[
\left(  \int \left( \int
(|DW(a,b,\theta)|a^{-s})^q
\frac{da}{a^{1+2q/p}}
\right)^{p/q}
  d\theta db \right)^{1/p}.
\]
In both cases we take $s = \sigma - 3/2(1/p-1/2)$.
(There is the possibility of defining
spaces using a third index (eg $B_{p,q,r}^\sigma$) corresponding
to the  $L^r$ norm in the $\theta$ variable,
but we ignore this here). As usual, the above formulas
can only provide norms for high-frequency functions,
and would have to be modified at coarse scales
if any low frequencies were present in $f$. As in the case of
the continuous wavelet transform for $\bf R$, there is some
heuristic value in considering the transform as measuring finite
directional differences e.g. $ f(b + a e_\theta) - 2 f(b) + f(b - a e_\theta)$,
where $e_\theta = (\cos(\theta),\sin(\theta))'$; however
this view is ultimately misleading. It is better to think
of the transform as comparing the difference between
polynomial approximation localized to two different
rectangles, one  of size $a$ by $\sqrt{a}$
and the other, concentric and co-oriented, of size
$2a$ by $\sqrt{2a}$.

The transform is
actually performing a kind of microlocal analysis of $f$
far more subtle than what is possible by simple
difference/differential expressions.
Indeed, consider the Heaviside $H(x) = 1_{\{x_1 > 0\}}$;
then at fine scales $DW(a,0,\theta) = 0$ for $|\theta| > \sqrt{a}$
and $DW(a,0,0) \approx a^{3/4}$ for $|\theta| \ll \sqrt{a}$, so that
$ \int_0^{2\pi} |DW(a,0,\theta)| d \theta \leq C a^{5/4}$
as $a \goto 0$. In short, $DW$ is giving very precisely the
orientation of the singularity.  Moreover, for $b \neq (0,x_2)'$,
$\int_0^{2\pi} |DW(a,0,\theta)| d \theta \goto 0$ rapidly as $a \goto 0$.
So the transform is localizing the singularity quite well
at fine scales, in a way that is difficult to imagine
simple differences
being able to do.
Interpreting the above observations, we learn that
a smoothly windowed Heaviside $f(x) = H(x) e^{- x^2}$
belongs in $DB_{\infty,\infty}^0$ but not
in any better space $DB_{\infty,\infty}^\sigma$, $\sigma > 0$,
while it belongs in $DB_{1,\infty}^1$ and not in any better space
$DB_{1,\infty}^\sigma$, $\sigma > 1$.  The difference between
the critical indices in these cases is indicative of the
sensitivity of the $p=1$ seminorms to sparsity.  Continuing in this
vein, we have that for weak $\ell^p$ embeddings, for each $\eta > 0$
\[
    \mu \{ (a,b,\theta) : |DW(a,b,\theta)| > \eps \} \leq C \eps^{-
(3/2 - \eta)}
\]
so that the space-scale-direction plane for the (windowed) Heaviside
is almost in $L^{2/3}(\mu) \sim DB_{2/3,2/3}^{3/2}$; the Heaviside
has something like $3/2$-derivatives.  In comparison, the wavelet expansion
of the Heaviside is only in $\ell^1$, so the expansion is denser and
`more irregular'
from the wavelet viewpoint than from the directional wavelet
viewpoint.  For comparison, the Dirac mass $\delta$
belongs at best to $B_{\infty,\infty}^{-1}$ and $B_{1,\infty}^{0}$
while it belongs at best to $DB_{\infty,\infty}^{-3/2}$ and
$DB_{1,\infty}^{-1/2}$.
The `point singularity' is more regular from the wavelet viewpoint than
from the directional wavelet viewpoint, while the Heaviside
is more regular from the directional wavelet viewpoint than
from the wavelet viewpoint. In effect, the Dirac `misbehaves in every
direction',
while the Heaviside misbehaves only in one direction, and this makes
a big difference
for the directional wavelet transform.

There are two obvious special equivalences:
first,  $L^2 \sim DF_{2,2}^0  \sim DB_{2,2}^0$
and $L^2$ Sobolev $W_2^m \sim DF_{2,2}^m  \sim DB_{2,2}^m$.
There are in general no other $L^p$ equivalences.
Outside the $L^2$ Sobolev scale, the only equivalence
with a previously
proposed space is with Hart Smith's ``Hardy Space for Fourier
Integral Operators'' \cite{SmithFIO}: ${\cal H}^1_{FIO} \sim DF_{1,2}^{0}$.
This space has a molecular decomposition
into {\it directional molecules}, which are functions that,
at high frequency, are
roughly localized
in space to an $a$ by $\sqrt{a}$ rectangle and roughly
localized in frequency to the dual rectangle rotated 90 degrees,
using traditional ways of measuring localization, such as
boundedness of moments of all orders in the two
principal directions.  Under this qualitative definition
of molecule, Smith showed that ${\cal H}^1_{FIO}$
has a molecular decomposition
$f =  \sum_Q A_Q m_Q(x)$ in which the coefficients
obey an $\ell^1$ norm summability condition $\sum |A_Q| 2^{j3/4} \leq 1$
when the directional molecules are $L^2$ normalized. This is obviously
the harbinger for a whole theory of directional molecular
decompositions.

More generally,
one can make a molecular decomposition of the directional
Besov and directional Triebel classes by discretizing
the directional wavelet transform according to
tiles $Q = Q(j,k_1,k_2,\ell)$ which obey the following desiderata:
\bitem
  \item In tile $Q(j,k_1,k_2,\ell)$,
   scale $a$ runs through a dyadic interval
   $ 2^{-j} > a \geq 2^{-(j+1)}$.
  \item At scale $2^{-j}$, locations run through rectangularly shaped
   regions with aspect ratio roughly $2^{-j}$ by $2^{-j/2}$.
  \item The location regions are rotated consistent with the orientation
\linebreak   $ b \approx R_{\theta_\ell} (k_1/2^{j},
k_2/2^{j/2})$.
  \item The tile contains orientations running through
  $ 2\pi \ell/2^{j/2} \leq \theta < 2\pi
  (\ell+1)/2^{j/2}$.
\eitem
Note again that for such tiles $\mu(Q) \approx 1$.
Over such tiles different values of $DW(a,b\theta)$ are roughly comparable
and different wavelets $\psi_{a,b,\theta}$ as well. Hence it is sensible to
decompose
\begin{eqnarray*}
    f(x) &=& \int DW(a,b,\theta) \psi_{a,b,\theta}(x) \mu(da db d\theta) \\
         &=& \sum_Q \int_Q DW(a,b,\theta) \psi_{a,b,\theta}(x) \mu(da
db d\theta)\\
         &=& \sum_Q M_Q(x), \qquad M_Q(x) = \int_Q DW(a,b,\theta)
\psi_{a,b,\theta}(x) \mu(da db d\theta) \\
         &=& \sum_Q A_Q m_Q(x), \qquad A_Q = \|DW(a,b,\theta) \|_{L^2(Q)}
\end{eqnarray*}

Morever, for any decomposition into directional molecules
(not just the approach above), the appropriate sequence norm of the
amplitude coefficients gives control of the corresponding
directional Besov or directional Triebel norm.
It is then relatively immediate that one can define
sequence space norms for which we have the
norm equivalences
\begin{equation}
     \| f \|_{DB_{p,q}^\sigma} \asymp  \| (A_Q)_Q
     \|_{db_{p,q}^\sigma},
\label{besovdirdef} \qquad
     \| f \|_{DF_{p,q}^\sigma} \asymp  \| (A_Q)_Q \|_{df_{p,q}^\sigma}
\end{equation}
where we again omit discussion of low frequency terms.
The sequence space equivalence $db_{2,2}^0 \sim df_{2,2}^0 \sim \ell^2$
are trivial.  An interesting equivalence
of relevance to the Heaviside example above is $ db_{2/3,2/3}^{3/2}
\sim \ell^{2/3}$,
so that, again, a smoothness space with
``$p < 1$'' is equivalent to an $\ell^\tau$ ball with $\tau < 1$.

Hart Smith made the crucial observation that
the molecules for the Smith space are
invariant under diffeomorphisms.
That is, if we take a $C^\infty$ diffeomorphism $\phi$,
and a family of ${\cal H}_{FIO}^1$
molecules (such as $m_Q(x)$),
then every $\tilde{m}_Q(x) = m_Q(\phi(x))$ is again a molecule, and
the sizes of moments defining the molecule property
are comparable for $m_Q$ and for $\tilde{m}_Q$.
It follows that ${\cal H}_{FIO}^1$ is invariant
under diffeomorphisms of the base space.
His basic lemma underlying this proof was strong enough
to apply to invariance of directional molecules in every
one of the directional Besov and directional Triebel classes.
Hence directional Besov and directional Triebel
classes are invariant under diffeomorphisms of the
base space.

This invariance enables a very simple calculation,
suggesting that the {\it directional wavelet transform sparsifies
objects with singularities along smooth curves}, or at least
sparsifies such objects to a greater extent that
does the ordinary wavelet transform. Suppose we analyse a function
$f$ which is smooth away from a discontinuity along a straight line;
then the Heaviside calculation we did earlier shows that
most directional wavelet coefficients are almost in weak $L^{2/3}$.
Now since objects with linear
singularities have $\ell^{2/3+\eps}$ boundedness of amplitudes
in a molecular decomposition, and directional molecules are
diffeomorphism invariant,
this sparsity condition
is invariant under diffeomorphisms
of the underlying space. It follows that
an object which is smooth away
from a discontinuity along
a smooth curve should also have molecular amplitudes
in $\ell^{2/3+\eps}$.

This sparsity argument suggests
that directional wavelets outperform wavelets for representing such
geometric objects.
Indeed, for $\eta > 0$ there is an $\eps > 0$
so that $\ell^{2/3+\eps}$ boundedness of directional
wavelet molecular amplitudes
shows that approximation by sums of $N$ directional molecules
allows a squared-$L^2$ approximation error of order $O(N^{-2+\eta})$,
whereas wavelet coefficients
of such objects are only in $\ell^1$, so
sums of $N$ wavelets only allow squared-$L^2$ approximation
error of size $O(N^{-1})$.

\subsection{Stylized applications}

\vskip-5mm \hspace{5mm}

The above calculations about sparsification of objects with curvilinear
singularities suggests the possibility of using the directional wavelet
transform based on parabolic scaling to pursue counterparts of all the
various classical wavelet applications mentioned in Section 3: nonlinear
approximation, data compression, noise removal, and fast computations.
It further suggests that such directional wavelet methods might
outperform calssical wavelets -- at least for objects containing
singularities along smooth curves, i.e. edges.

\subsubsection{First discretization: curvelets}

\vskip-5mm \hspace{5mm}

To develop applications, molecular decomposition is (once again) not
enough: some sort of rigid decomposition needs to be developed; an
orthobasis, for example.

Cand\`es and Donoho \cite{Curvelets} developed a
tight frame of elements exhibiting parabolic dilations
which they called {\it curvelets}, and used it to
systematically develop some of these applications.
A side benefit of their work is knowledge that the
transform is essentially optimal, i.e. that there
is no fundamentally better scheme of nonlinear approximation.
The curvelet system has a countable collection of generating
elements
$\gamma_\mu(x_1,x_2)$ , $\mu \sim (a_j, b_{k_1,k_2}, \theta_{\ell}, t_{m})$
which code for scale, location, and direction.
They obey the usual rules for a tight frame, namely,
the reconstruction formula
and the Parseval relation:
\[
    f = \sum_\mu \langle \gamma_\mu, f \rangle \gamma_\mu ,  \qquad
    \| f \|_2^2 = \sum_\mu \langle \gamma_\mu, f \rangle^2.
\]
The transform is based on a series of space/frequency
localizations, as follows.
   \bitem
     \item Bandpass filtering.  The object is separated
     out into different dyadic scale subbands, using
     traditional bandpass filtering with passband
     centered around $|\xi| \in [2^j, 2^{j+1}]$.
     \item Spatial localization.  Each bandpass object
     is then smoothly partitioned spatially into boxes of side
       $2^{-j/2}$.
     \item Angular localization.  Each box is analysed
     by ridgelet transform.
   \eitem

The frame elements are
essentially localized into boxes of side $2^{-j}$ by $2^{-j/2}$
at a range of scales, locations, and orientations,
so that it is completely consistent with the molecular
decomposition of the directional Besov or directional
Fourier classes.  However, unlike the molecular decomposition,
the coefficients are
linear in $f$ and the frame elements are fixed elements.
Moreover, an algorithm for application to real data
on a grid is relatively immediate.

\subsubsection{Nonlinear approximation}

\vskip-5mm \hspace{5mm}

In dimension 2, the analog to what was called free knot spline
approximation is approximation by piecewise polynomials on
triangulations with $N$ pieces. This idea has generated a lot of
interest but frustratingly few hard results. For one thing, it is not
obvious how to build such triangulations in a way that will fulfill
their apparent promise, and in which the resulting algorithm is
practical and possible to analyze.

Here is a class of two-dimensional
functions where this scheme
might be very attractive.
Consider a class $\cal F$ of model `images' which exhibit
discontinuities across
$C^2$ smooth curves. These `images'
are supposed to be $C^2$ away from discontinuity.
Moreover, we assume uniform control
both of the $C^2$ norm for the discontinuity curve and
smooth function.  One can imagine that very fine
needle-like triangles near curved discontinuities
would be valuable; and this is indeed so, as \cite{SCA}
shows; in an ideal triangulation one geta a squared error converging at rate
$N^{-2}$ whereas adaptive quadtrees and other simpler partitioning
schemes give only $N^{-1}$ convergence.
Moreover, this rate is optimal, as shown in \cite{SCA},
if we allow piecewise smooth approximation on
essentially arbitrary triangulations with $N$ pieces,
even those designed
by some as yet unknown very clever and very nonlinear
algorithm, we cannot in general converge to such objects
faster than rate $N^{-2}$.

Surprisingly, a very concrete algorithm does almost this well: simply
thresholding the curvelet coefficients. Cand\`es and Donoho have shown
the following \cite{CurveletMSS}

{\bf Theorem}: {\it The decreasing rearrangement of the frame
coefficients in the curvelet system obeys the following inequality for
all $f \in \cal F$:}
\[
    |\alpha|_{(k)} \leq C k^{-3/2} \log^{3/2} (k), \qquad k \geq 1.
\]

This has exactly the implication one would have
hoped for from the molecular decomposition of
directional Besov classes: the frame coefficients
are in $\ell^{2/3+\eps}$ for each $\eps > 0$.
Hence, we can build an approximation to
a smooth object with
curvilinear discontinuity from $N$ curvelets with squared $L^2$-error
$\log^{3}(N) \cdot N^{-2}$; as mentioned earlier,
Wavelets would give squared $L^2$-error $\geq c N^{-1}$.

In words: approximation by sums of the $N$-biggest
curvelet terms does essentially
as well in approximating objects in $\cF$
as free-triangulation into $N$ regions. In a sense, the
result is analogous to the result mentioned
above in Section 3.2.1 comparing wavelet thresholding
to nonlinear spline approximation,
where we saw that approximation by
the $N$-biggest amplitude wavelet terms does
as well as free-knot splines with $N$ knots.
There has been a certain amount of talk about
the problem of characterizing approximation spaces for
approximation by $N$ arbitrary triangles; while
this problem seems very intractable, it is clear
that the directional Besov classes provide
what is, at the moment, the next best thing.

\subsubsection{Data compression}

\vskip-5mm \hspace{5mm}

Applying just the arguments already given in the wavelet case show that
the result of $L^2$ nonlinear approximation by curvelets, combined with
simple quantization, gives near-optimal compression of functions in the
class $\cal F$ above, i.e. the number of bits in the compressed
representation is optimal to within some polylog factor. This seems to
promise some interesting practical coders someday.

\subsubsection{Noise removal}

\vskip-5mm \hspace{5mm}

The results on nonlinear approximation by thresholding of the curvelet
coefficients have corresponding implications in statistical estimation.
Suppose that we have noisy data according to the white noise model
\[
   Y(dx_1,dx_2) = f(x_1,x_2) dx_1 dx_2 + \eps W(dx_1,dx_2)
\]
where $W$ is a Wiener sheet. Here $f$ comes from the same `Image Model'
$\cF$ discussed earlier, of smooth objects with discontinuities across
$C^2$ smooth curves. We measure risk by Mean Squared Error, and consider
the estimator that thresholds the curvelet coefficients at an
appropriate (roughly $2\sqrt{\log(\eps^{-1})}$) multiple of the noise
level. Emmanuel Cand\`es and I showed the following
\cite{CurveletsRadon}:

{\bf Theorem}: {\it Appropriate thresholding of curvelet coefficients
gives nearly the optimal rate of convergence; with $polylog(\eps)$ a
polynomial in $\log(1/\eps)$, the estimator $\hat{f}^{CT}$ obeys}
\[
    R_\eps( f, \hat{f}^{CT} ) \leq polylog(\eps) \cdot \min_{\hat{f}} \max_{f
\in \cF} R_\eps (f,\hat{f}).
\]

Hence, in this situation, curvelet thresholding
outperforms wavelet thresholding at the level of rates:
$O(polylog(\eps) \cdot \eps^{4/3})$ vs $O(\eps)$.
Similar results can be developed for other
estimation problems, such as the problem of Radon
inversion. There the rate comparison is
$polylog(\eps) \cdot \eps^{4/5}$ vs
$\log(1/\eps) \cdot \eps^{2/3}$; \cite{CurveletsRadon}.
In empirical work \cite{StarckIEEE,CandesGuo}, we have seen
visually persuasive results.

\subsubsection{Improved discretization: directional framelets}

\vskip-5mm \hspace{5mm}

The curvelet representation described earlier is a somewhat awkward way
of obtaining parabolic scaling, and also only indirectly related to the
continuum directional wavelet transform. Cand\`es and Guo
\cite{CandesGuo} suggested a different tight frame expansion based on
parabolic scaling.  Although this was not introduced in such a fashion,
for this exposition, we propose an alternate way to understand their
frame, simply as discretizing the directional wavelet transform in a way
reminiscent of (\ref{orthodef}); for details, see \cite{DWTURL}.
Assuming a very specific choice of directional wavelet, one can get (the
fine scale) frame coefficients simply by sampling the directional
wavelet transform, obtaining a decomposition
\[
      f = \sum_{j,k,l} DW(2^{-j},b_{k_1,k_2}^{j,\ell},2\pi l/2^{j/2})
\psi_{2^{-j},k/2^j,2\pi\ell/2^{j/2}} =
\sum_{j,k,l}
\alpha_{j,k,l}
\psi_{j,k,l}, \mbox{say } ;
\]
(as usual, this is valid as written only for high-frequency functions).
In fact this can yield a tight frame, in particular the Parseval relation
$\sum_{j,k,l} \alpha_{j,k,l}^2 = \| f \|_{L^2}^2 $.
This has conceptual advantages: a better relationship
to the continuous directional wavelet transform and
perhaps an easier path to digital representation.
In comparison with the original curvelets scheme,
curvelets most naturally organizes matters so that `within' each location
we see all directional behavior represented, whereas directional
framelets most naturally organize matters so that `within'
each orientation we see all locations represented.

\subsubsection{Operator representation}

\vskip-5mm \hspace{5mm}

Hart Smith, at the Berlin ICM, mentioned that decompositions based on
parabolic scaling were valuable for understanding Fourier Integral
Operators (FIO's) \cite{SmithICM}; in the notation of our paper, his
claim was essentially that FIO's of order zero operate on fine-scale
directional molecules approximately by performing well-behaved affine
motions -- roughly, displacement, scaling and change of orientation.
Underlying his argument was the study of families of elements generated
from a single wavelet by true affine parabolic scaling
$\phi_{a,b,\theta}(x) = \phi( P_a  \circ  R_\theta \circ S_b x )$ where
$P_a = diag(a,\sqrt{a})$ is the parabolic scaling operator and $S_b x =
x-b $ is the shift. Smith showed that  if $T$ is an FIO of order $0$ and
$\phi$ is directionally localized, the kernel
\[
     K_{a,b,\theta}^{a',b',\theta'} = \langle \phi_{a,b,\theta} , T
\phi_{a',b',\theta'} \rangle
\]
is rapidly decaying in its entries as one moves away from `the
diagonal' in an appropriate sense.

Making this principle more adapted to discrete frame representations
seems an important priority.
Cand\`es and Demanet have recently
announced \cite{CandesDemanet}
that actually, {\it
the matrix representation of FIOs of order $0$ in
the directional framelet decomposition is
sparse.} That is, each row and column of the
matrix will be in $\ell^p$ for each $p > 0$.
in a directional wavelet frame.
This observation is analogous in some
ways to Meyer's observation that the
orthogonal wavelet transform
gives a sparse representation for Calder\'on-Zygmund operators.
Cand\`es has hopes that this sparsity may form some day the
basis for fast algorithms for hyperbolic PDE's and other FIO's.

\subsection{Applications}

\vskip-5mm \hspace{5mm}

The formalization of the directional wavelet transform and curvelet
transform are simply too recent to have had any substantial applications
of the `in daily use by thousands' category.  Serious deployment into
applications in data compression or statistical estimation is still off
in the future.

However, the article
\cite{DonohoFlesia}  points to the possibility of immediate effects on research
activity in computational neuroscience, simply
by generating new research hypothesis. In effect,
if vision scientists can be induced
to consider these new types of image representation,
this will stimulate
meaningful new experiments, and re-analyses of existing
experiments.

To begin with, for decades, vision scientists
have been influenced by
mathematical ideas in framing research hypotheses
about the functioning of the visual cortex,
particular the functioning of the V1 region.
In the 1970's, several authors suggested
that the early visual
system does Fourier Analysis;
by the 1980's the cutting edge hypothesis became the
suggestion that the
early visual system does
Gabor Analysis; and by the 1990's, one
saw claims that the
early visual system does
a form of wavelet analysis.
While the hypotheses have changed over time,
the invariant is that
vision scientists have relied on mathematics to
provide language \& intellectual framework
for their investigations.
But it seems likely that the hypotheses of these previous
decades are incomplete, and that to these should be
added the hypothesis that the early visual system
performs a directional wavelet transform
based on parabolic scaling.  During my Plenary Lecture,
biological evidence was presented consistent
with this hypothesis, and a proposal was made that
future experiments in intrinsic optical imaging
of the visual cortex ought to attempt to test
this hypothesis. See also \cite{DonohoFlesia}.

\section{Geometric multiscale analysis `without Calder\'on'}

\vskip-5mm \hspace{5mm}

In the last section we considered a kind of geometric
multiscale analysis employing a Calder\'on-like formula.
In the ICM Lecture we also considered
dispensing with the need for Calder\'on formulas,
using a cruder set of multiscale tools,
but one which allows for a wide range of interesting
applications -- very different from the applications
based on analysis/synthesis and Parseval. Our model
for how to get started in this direction
was Peter Jones' travelling salesman
problem. Jones considered instead a
countable number of points $X = \{x_i\}$ in $[0,1]^2$
and asked: when can the points of $X$ be connected
by a finite length (rectifiable) curve ?
And, if they can be,
what is the shortest possible length?
Jones showed that  one should consider,
for each dyadic square
$Q$ such that the dilate $3Q$ intersects $X$, the width $w_Q$
of the thinnest strip in the plane containing
all the points in $X \cap 3Q$, and define $\beta_Q = w_Q/diam(Q)$
the proportional width of that strip, relative to the sidelength
of $Q$. As $\beta_Q = 0$ when the data lie on a straight line,
this is precisely a
measure of how close to linear
the data are over the square $Q$. He proved the
there is a finite-length curve $\Gamma$ visiting all
the points in $X = \{x_i\}$ iff
$\sum \beta_Q^2 diam(Q)  < \infty $.  I find it very impressive
that analysis of the number of points in strips of various
widths can reveal the existence of a rectifiable curve
connecting those points. In our lecture,
we discussed this idea of counting points in anistropic strips
and several applications in signal detection and pattern recognition
\cite{ConnectDots,MultiscaleDetect},
with applications in characterizing galaxy clustering \cite{Galaxies}.
We also referred to interesting work such as Gilad Lerman's
thesis \cite{LermanThesis}, under the direction of Coifman and
Jones, and to \cite{BMIA}, which surveys a wide range of
related work.  Look to \cite{DWTURL} for an extended
version of this article covering such topics.

\section{Conclusion}

\vskip-5mm \hspace{5mm}

Important developments in `pure' harmonic analysis,
like the use of parabolic scaling for
study of convolution operators and FIOs,
or the use of anisotropic strips
for analysis of rectifiable measures, did
not arise because of applications
to our developing `information society', yet they seem to
have important stylized applications
which point clearly in that direction.  A number
of enthusiastic applied mathematicians,
statisticians, and scientists
are attempting to develop true `real world' applications.

At the same time, the fruitful directions
for new kinds of geometric multiscale analysis
and the possible limitations to be surmounted
remain to be determined.
Stay tuned!

\section{Acknowledgements}

\vskip-5mm \hspace{5mm}

The author would like to thank Emmanuel Cand\`es, Raphy Coifman, Peter
Jones, and Yves Meyer for very insightful and inspiring discussions.
This work has been partially supported by National Science Foundation
grants DMS 00-77261, 98--72890 (KDI), and DMS 95--05151.

\label{lastpage}

\end{document}